\documentclass[12pt]{amsart}
\usepackage{amssymb}
\usepackage{amsfonts}
\usepackage{amssymb,latexsym}
\usepackage{enumerate}
\usepackage{url}
\usepackage{mathrsfs}
\makeatletter
\@namedef{subjclassname@2010}{%
  \textup{2010} Mathematics Subject Classification}
\makeatother

  [2002/01/19 v2.2g %
    AMS font definitions%
  ]
\DeclareFontFamily{U}{euf}{}
\DeclareFontShape{U}{euf}{m}{n}{%
  <5><6><7><8><9>gen*eufm%
  <10><10.95><12><14.4><17.28><20.74><24.88>eufm10%
  }{}
\DeclareFontShape{U}{euf}{b}{n}{%
  <5><6><7><8><9>gen*eufb%
  <10><10.95><12><14.4><17.28><20.74><24.88>eufb10%
  }{}

  [2002/01/19 v2.2g %
    AMS font definitions%
  ]
\DeclareFontFamily{U}{msb}{}
\DeclareFontShape{U}{msb}{m}{n}{%
  <5><6><7><8><9>gen*msbm%
  <10><10.95><12><14.4><17.28><20.74><24.88>msbm10%
  }{}

  [2002/01/19 v2.2g %
    AMS font definitions%
  ]
\DeclareFontFamily{U}{msa}{}
\DeclareFontShape{U}{msa}{m}{n}{%
  <5><6><7><8><9>gen*msam%
  <10><10.95><12><14.4><17.28><20.74><24.88>msam10%
  }{}

\newtheorem{theorem}{Theorem}[section]

\newtheorem{proposition}[theorem]{Proposition}
\newtheorem{corollary}[theorem]{Corollary}

\theoremstyle{definition}

\newtheorem{remark}[theorem]{Remark}

\numberwithin{equation}{section} 

\begin{document}

\title[Generalizations of Lerch's formula]{Generalizations of Lerch's formula by Barnes' multiple zeta functions}

\author{Su Hu}
\address{Department of Mathematics, South China University of Technology, Guangzhou, Guangdong 510640, China}
\email{mahusu@scut.edu.cn}

\author{Min-Soo Kim}
\address{Department of Mathematics Education, Kyungnam University, Changwon, Gyeongnam 51767, Republic of Korea}
\email{mskim@kyungnam.ac.kr}



\begin{abstract}
The classical Lerch's formula states the following normalized product:
$$\prod_{n=0}^\infty(x+n)=\frac{\sqrt{2\pi}}{\Gamma(x)},\quad \textrm{Re}(x)>0,
$$
where $\Gamma(x)$ is the Euler gamma function.

In this note, by using Barnes' multiple zeta function and its alternating form, 
we obtain two kinds of generalizations of Lerch's formula,  which imply the product
$$
\prod_{n=1}^\infty n=\sqrt{2\pi}
$$
(in the sense of zeta regularization)
and the product
$$\frac{2\cdot2}{1\cdot 3}\frac{4\cdot4}{3\cdot 5}\frac{6\cdot6}{5\cdot 7}\cdots=\frac{\pi}{2}$$
(Wallis' formula in 1656), respectively.

\end{abstract}

\maketitle

\def\ord{\text{ord}_p}
\def\ordt{\text{ord}_2}
\def\o{\omega}
\def\la{\langle}
\def\ra{\rangle}
\def\Log{{\rm Log}\, \Gamma_{p,N}^*}
\def\ov{\bar}

\section{Introduction}
In the book ``Number theory, 3, Iwasawa theory and modular forms" \cite{KKS}, the authors stated the following infinite product from the zeta regularization:
\begin{equation}\label{Riemann}
\infty !:=\prod_{n=1}^\infty n=\sqrt{2\pi}.
\end{equation}
(See \cite[Corollary 9.13]{KKS}).

They explained it from the   following interesting way.
Given a finite sequence $\textbf{a}=(a_{1},a_{2},\ldots, a_{N}),$
define  $$\zeta_{\textbf{a}}(s)=\sum_{n=1}^{N}a_{n}^{-s}.$$
Then $$\zeta_{\textbf{a}}^{'}(0)=-\sum_{n=1}^{N}\log (a_{n})=-\log\left(\prod_{n=1}^{N}a_N\right).$$
So $$\exp(-\zeta_{\textbf{a}}^{'}(0))=\prod_{n=1}^{N}a_{n}.$$
Thus for an infinite sequence $\textbf{a}=(a_{1},a_{2},\ldots, a_{N},\ldots),$
if the corresponding zeta function 
$$\zeta_{\textbf{a}}(s)=\sum_{n=1}^{\infty}a_{n}^{-s}$$
can be analytically continued to a neighborhood around $s=0$,
then we may define their normalized product or zeta regularization by
\begin{equation}\label{normal}
\prod_{n=1}^{\infty}a_{n}=\exp(-\zeta_{\textbf{a}}^{'}(0)).
\end{equation}
Now let  $\textbf{a}=(1,2,3,\ldots),$
then the corresponding zeta function is just the Riemann zeta function
$$\zeta(s)=\sum_{n=1}^{\infty}\frac{1}{n^{s}},$$
and in 1859, Riemann \cite{Riemann} proved that 
\begin{equation}\label{Riemann2}
\zeta^{'}(0)=-\frac{1}{2}\log(2\pi).
\end{equation}
By (\ref{normal}), we understand that
\begin{equation}
\prod_{n=1}^{\infty}n=\exp(-\zeta^{'}(0))=\sqrt{2\pi}.
\end{equation}

In 1894, Lerch generalized (\ref{Riemann}) as the following product:
\begin{equation}\label{Lerch}
\prod_{n=0}^\infty(x+n)=\frac{\sqrt{2\pi}}{\Gamma(x)},\quad \textrm{Re}(x)>0,
\end{equation}
where $\Gamma(x)$ is the Euler gamma function.
(See {\cite[Theorem 9.12]{KKS}). In fact, letting $x=1$ in (\ref{Lerch}),
we get (\ref{Riemann}).

Denote by $$\mathbb{Z}_{0}^{-}=\{0,-1,-2, \ldots\}.$$
In 2004, the above Lerch's formula was generalized by Kurokawa and Wakayama \cite{KW} in the following way:
\begin{equation}
\prod_{n=0}^{\infty}\left((n+x)^{m}-y^m\right)=\frac{(\sqrt{2\pi})^{m}}{\prod_{\zeta^{m}=1}\Gamma(x-\zeta y)}
\end{equation}
and in 2006, applying Stark's summation formula \cite{Stark}, Mizuno \cite{Mizuno} got  a further extension:
\begin{equation}
\prod_{m=0}^{\infty}\left(\prod_{j=1}^{n}(m+z_{j})\right)=\frac{(\sqrt{2\pi})^{n}}{\prod_{j=1}^{n}\Gamma(z_{j})}=\prod_{j=1}^{n}\left(\prod_{m=0}^{\infty}(m+z_{j})\right)
\end{equation}
for $z_{j}\in\mathbb{C}\setminus\mathbb{Z}_{0}^{-}.$

In this note,  
we obtain two kinds of generalizations of Lerch's formula by using Barnes' multiple zeta function and its alternating form, respectively.

Suppose that $\omega_1,\ldots,\omega_N$ are positive real numbers and $x$ is a complex number with positive real part. 
In 1904, Barnes \cite{Ba} studied the multiple Hurwitz zeta function  
\begin{equation}\label{Barnes}
\zeta_{N}(s,x;\omega_1,\ldots,\omega_N)=\sum_{m_1=0}^\infty\cdots\sum_{m_N=0}^\infty
\frac{1}{(x+\omega_1m_1+\cdots+\omega_Nm_N)^s} ,\quad\text{Re}(s)>N,
\end{equation}
and the corresponding gamma function $\Gamma_{N}^{B}(x).$

Setting $$\omega_{j}=1,~~ \textrm{for}~~ j=1,2,3,\ldots,N$$ in (\ref{Barnes}),
we get 
the following special case:  
\begin{equation}\label{M-Z}
\zeta_{N}(s,x)=\sum_{m_1,\ldots,m_N=0}^\infty
\frac{1}{(x+m_1+\cdots+m_N)^s},\quad\text{Re}(s)>N, x\in\mathbb{C}\setminus\mathbb{Z}_{0}^{-}.
\end{equation} 
As shown by Choi and Srivastava in \cite[p. 503--504]{Choi} or  \cite[p. 149, Theorem 2.6]{SC}, for Re$(x)>0,$ the function $\zeta_{N}(s,x)$ can be analytic continued  for all $s\in\mathbb{C}$ except for simple poles at
$$s=k  \quad(1\leq k\leq N; N,k\in\mathbb{N}).$$
Thus the multiple (or, simply, $N$-ple) gamma function $\Gamma_{N}(x)$ can now be defined by
\begin{equation}\label{ru1}
\Gamma_{N}(x)=\exp\left(\frac{\partial}{\partial s}\zeta_{N}(s,x)\biggl|_{s=0}\right), \quad\text{Re}(x)>0.
\end{equation}
(See \cite{KK}, \cite{Ru} and \cite{SC}).
The following recurrence formula of $\Gamma_{N}(x)$ is well-known:
\begin{equation}\label{recurrence}
\Gamma_{N}(x+1)=\frac{\Gamma_{N}(x)}{\Gamma_{N-1}(x)}
\end{equation} 
(e.g. see \cite[Proposition 1.2 (1)]{Choi3})
and several series expansions and asymptotic formulas of $\Gamma_{N}(x)$ and $\log\Gamma_{N}(x)$ have been considered by Choi and Srivastava in \cite[p. 375-382]{SC}.

By using  the multiple zeta function $\zeta_{N}(s,x)$ and the corresponding gamma function  $\Gamma_{N}(x)$,  we first prove the following generalization
of the above Lerch's formula (\ref{Lerch}).
\begin{theorem}[The first generalized Lerch formula]\label{gen-lerch}
For $N\in\mathbb{N},$ in the sense of zeta regularization, we have
\begin{equation}\label{GLerch}
\prod_{n=0}^\infty(x+n)^{-\binom{n+N-1}{N-1}}=\Gamma_{N}(x), \quad {\rm Re}(x)>0.\end{equation}
\end{theorem}
  The following relation between  $\Gamma_{1}(x)$ and  the Euler gamma function $\Gamma(x)$ will be proved in Section \ref{B}. 
 \begin{proposition}[{\cite[Lemma 2.1]{Var}}]\label{proposition1}
 \begin{equation}\label{G1-re}
\Gamma_1(x)=\frac{\Gamma(x)}{\sqrt{2\pi}}.
\end{equation}
\end{proposition} 
 From the above proposition, we see that if letting $N=1$ in (\ref{GLerch}), then we recover Lerch's formula (\ref{Lerch}).
 
Now we go to the alternating case. 
For  $\textrm{Re}(s)>0$, the multiple Barnes-Euler zeta function $\zeta_{E,N}(s,x;\omega_1,\ldots,\omega_N)$   is
defined as a deformation of  the Barnes' multiple  zeta function as follows
\begin{equation}\label{Barnes2}
\zeta_{E,N}(s,x;\omega_1,\ldots,\omega_N)=\sum_{m_1=0}^\infty\cdots\sum_{m_N=0}^\infty
\frac{(-1)^{m_1+\cdots+m_N}}{(x+\omega_1m_1+\cdots+\omega_Nm_N)^s}.
\end{equation}
Its analytic properties both in the complex plane $\mathbb{C}$ and in the $p$-adic complex plane $\mathbb{C}_{p}$ have been systematically 
studied by the authors in \cite{arXiv}.

Letting $$\omega_{j}=1,~~ \textrm{for}~~ j=1,2,3,\ldots,N$$ in (\ref{Barnes2}),
we get 
the following special case:  
\begin{equation}\label{M-Z-2}
\zeta_{E, N}(s,x)=\sum_{m_1=0}^\infty\cdots\sum_{m_N=0}^\infty
\frac{(-1)^{m_1+\cdots+m_N}}{(x+m_1+\cdots+m_N)^s},\quad\text{Re}(s)>0, x\in\mathbb{C}\setminus\mathbb{Z}_{0}^{-}.
\end{equation} 
As shown by Choi and Srivastava in \cite[Section 3]{Choi}, for Re$(x)>0,$ the function $\zeta_{E, N}(s,x)$ can be continued analytically as an entire function of $s\in\mathbb{C}.$  
Thus the corresponding   multiple (or, simply, $N$-ple) gamma function $\Gamma_{N}^{*}(x)$ can now be defined by
\begin{equation}\label{ru2}
\Gamma_{N}^*(x)=\exp\left(\frac{\partial}{\partial s}\zeta_{E,N}(s,x)\biggl|_{s=0}\right), \quad\text{Re}(x)>0.\end{equation}
The following recurrence formula of $\Gamma_{N}^{*}(x)$ is implied by a formula of the authors~\cite[Lemma 2.1 (1)]{arXiv}:
\begin{equation} 
\Gamma_{N}^{*}(x+1)=\frac{\Gamma_{N-1}^{*}(x)}{\Gamma_{N}^{*}(x)}.
\end{equation} 
By using  the alternating multiple zeta function $\zeta_{E,N}(s,x)$ and the corresponding gamma function  $\Gamma_{N}^{*}(x)$,  we prove another generalization
of the above Lerch's formula (\ref{Lerch}).
\begin{theorem}[The second gneralized Lerch formula]\label{gen-lerch2}
For $N\in\mathbb{N},$ we have
\begin{equation}\label{GLerch2} 
\prod_{n=0}^\infty(x+n)^{(-1)^{n+1}\binom{n+N-1}{N-1}}=\Gamma_{N}^*(x),\quad {\rm Re}(x)>0.\end{equation}
\end{theorem}
The following relation between  $\Gamma_{1}^{*}(x)$ and  the Euler gamma function $\Gamma(x)$ will be proved in Section \ref{EB}.
\begin{proposition}\label{proposition2}
\begin{equation}\label{ex-gamma}
\Gamma_{1}^*(x)=\frac{\Gamma\left(\frac{x}{2}\right)}{\sqrt{2}\Gamma\left(\frac{x+1}{2}\right)},\quad{\rm Re}(x)>0.
\end{equation}
\end{proposition}
\begin{remark}
In \cite[p. 742]{Miller}, during the proof for the claim that the Students $t$-distribution is a continuous probability density, Miller obtained the following integral
\begin{equation}\label{Miller}
\int_{-\infty}^{\infty} \left(1+\frac{t^{2}}{x}\right)^{-\frac{x+1}{2}} dt=\frac{\sqrt{\pi x}\Gamma\left(\frac{x}{2}\right)}{\Gamma\left(\frac{x+1}{2}\right)}.
\end{equation}
From Proposition \ref{proposition2}, the above equality in fact leads to the following integral representation of $\Gamma_{1}^*(x)$:
\begin{equation}\label{Miller2}
\Gamma_{1}^*(x)=\frac{\Gamma\left(\frac{x}{2}\right)}{\sqrt{2}\Gamma\left(\frac{x+1}{2}\right)}=\frac{1}{\sqrt{2\pi x}}\int_{-\infty}^{\infty}\left(1+\frac{t^{2}}{x}\right)^{-\frac{x+1}{2}} dt.
\end{equation}
\end{remark}
By setting $N=1$ in (\ref{GLerch2}), from (\ref{ex-gamma}), we immediately get
\begin{corollary}[Lerch type formula]
\begin{equation}\label{Lerch2}
\prod_{n=0}^\infty(x+n)^{(-1)^{n+1}}=\frac1{\sqrt2}\frac{\Gamma\left(\frac{x}{2}\right)}{\Gamma\left(\frac{x+1}{2}\right)},\quad {\rm Re}(x)>0.
\end{equation}
\end{corollary}
Then letting $x=1$ in (\ref{Lerch2}), and recalling that $\Gamma(1)=1$ and $\Gamma(\frac{1}{2})=\sqrt{\pi}$, we immediately get
\begin{corollary}[Wallis' formula]
\begin{equation}
\frac{2\cdot2}{1\cdot 3}\frac{4\cdot4}{3\cdot 5}\frac{6\cdot6}{5\cdot 7}\cdots=\frac{\pi}{2}.
\end{equation}
\end{corollary}
\begin{remark} As early as in 1656, the British mathematician John Wallis \cite{Wallis} showed the above remarkable formula  in his book  ``Arithmetica Infinitorum".
 After his work, there appears  many methods to prove it, including the well-known ones based on the formula for integrals of powers of $\sin x$ from the inductive method 
  or based on  the infinite product expansion of $\sin x$. Recently, Miller \cite{Miller} found a probabilistic proof by using the Students $t$-distribution, 
and Friedmann and Hagen \cite{Friedmann} presented an quantum mechanical derivation based on  the spectrum of the hydrogen in the physical three dimensions.
In 1994, Sondow \cite{Sondow} derived (\ref{Riemann2})  from Wallis' formula by using Euler's  transformation of series.
\end{remark}

\section{Barnes' multiple zeta function and the first generalized Lerch  formula}\label{B}
The main aim of this section is to prove the first generalized Lerch's formula (Theorem \ref{gen-lerch} above).

\subsection*{Proof of Theorem \ref{gen-lerch}}

According to Barnes \cite{Ba}, as the classical Riemann zeta functions $\zeta(s),$ the multiple zeta function $\zeta_{N}(s,x)$ defined in (\ref{M-Z}) may also be represented by the Mellin transform as follows
\begin{equation}\label{B-E-def-3}
\zeta_{N}(s,x)=\frac1{\Gamma(s)}\int_0^\infty t^se^{-xt}\left(1-e^{-t}\right)^{-N}\frac{dt}{t},
\end{equation}
for Re$(s)>N$ and Re$(x)>0.$

This can be shown as follows.
Start with the  power series 
\begin{equation}\label{re-1}
\begin{aligned}
\left(1-e^{-t}\right)^{-N}=\sum_{m_1,\ldots,m_N=0}^\infty e^{-t(m_1+\cdots+m_N)},\quad t>0.
\end{aligned}
\end{equation}
Note that 
\begin{equation}\label{re-2}
w^{-s}=\frac1{\Gamma(s)}\int_0^\infty t^s e^{-w t}\frac{dt}{t} 
=\frac1{\Gamma(s)}\int_1^\infty u^{-w}(\log u)^{s-1}\frac{du}{u}
\end{equation}
(see \cite[(2.7)]{Ru} and \cite[(1)]{SC}).
Then substituting  (\ref{re-1}) into (\ref{B-E-def-3}) we have
\begin{equation}\label{B-E-def-2}
\begin{aligned}
\zeta_{N}(s,x)
&=\frac1{\Gamma(s)}\int_0^\infty t^se^{-xt}\left(1-e^{-t}\right)^{-N}\frac{dt}{t}
\\&=\sum_{m_1,\ldots,m_N=0}^\infty\frac1{\Gamma(s)}\int_0^\infty t^s e^{-(x+m_1+\cdots+m_N) t}\frac{dt}{t} \\
&=\sum_{m_1,\ldots,m_N=0}^\infty\frac{1}{(x+m_1+\cdots+m_N)^s},
\end{aligned}
\end{equation}
for Re$(s)>N$ and Re$(x)>0.$ (See \cite[(3.14)]{Ru}).

By (\ref{B-E-def-3}) we have
\begin{equation}\label{dir-mul-int}
\begin{aligned}
\zeta_{N}(s,x)&=\frac1{\Gamma(s)}\int_0^\infty t^se^{-xt}\left(1-e^{-t}\right)^{-N}\frac{dt}{t} \\
&=\frac1{\Gamma(s)}\int_1^\infty\left(\frac1{1-u^{-1}}\right)^{N}u^{-x}(\log u)^{s-1}\frac{du}{u} \\
&=\frac1{\Gamma(s)}\sum_{n=0}^\infty(-1)^n\binom{-N}{n}\int_1^\infty u^{-x-n}(\log u)^{s-1}\frac{du}{u} \\
&\quad(\text{by using (\ref{re-2}) with $w=x+n$}) \\
&=\sum_{n=0}^\infty \binom{n+N-1}{N-1}(x+n)^{-s}
\end{aligned}
\end{equation}
(see \cite[p. 505, (1.17)]{Choi}). The above equation is established for Re$(x)>0$ and  for all $s\in\mathbb{C}$ except  
$$s=k  \quad(1\leq k\leq N; N,k\in\mathbb{N}).$$ In fact, (\ref{dir-mul-int}) is another way of analytic continuation for $\zeta_{N}(s,x)$ by Choi in \cite{Choi2}. 
Then from (\ref{ru1}), in considering of the concept of normalized product or zeta regularization introduced at the beginning of this paper, we have
\begin{equation}\label{dir-mul-int-an}
\begin{aligned}
\Gamma_{N}(x)=&\exp\left(\frac{\partial}{\partial s}\zeta_{N}(s,x)\biggl|_{s=0}\right)
\\&=\exp\left(-\sum_{n=0}^\infty\binom{n+N-1}{N-1}\log(x+n)\right) \\
&=\prod_{n=0}^\infty(x+n)^{-\binom{n+N-1}{N-1}},
\end{aligned}
\end{equation}
which is the desired result.

\subsection*{Proof of Proposition \ref{proposition1}}

In the history, the gamma function $\Gamma(s)$ is defined by Euler in 1729 from the integral
\begin{equation}\label{gamma}
\Gamma(s)=\int_0^\infty e^{-t}t^{s} \frac{dt}{t}.
\end{equation}
Note that this integral is well-defined if Re$(s)>0.$  Integrating by parts yields $$\Gamma(s+1)=s\Gamma(s).$$
This implies that if $n$ is a nonnegative integer then $$\Gamma(n+1)=n!,$$ thus $\Gamma(s)$ generalizes the factorial function.

Letting $N=1$ in (\ref{M-Z}), we recover the Hurwitz zeta function introduced by Hurwitz in 1882:
$$\zeta(s,x)=\sum_{n=0}^{\infty}\frac{1}{(n+x)^{s}}.$$

By (\ref{ru1}), we have \begin{equation}\label{Ga}\Gamma_1(x)=e^{\zeta'(0,x)}.\end{equation}
The following equality comes from the  difference functional equation of the Hurwitz zeta function $\zeta(s,x)$: 
\begin{equation}\label{zeta}\zeta'(0,x+1)=\zeta'(0,x)+\log x.\end{equation} (See \cite[p. 497]{Var}). 
Then substituting (\ref{zeta}) into (\ref{Ga}), we get 
\begin{equation}\label{G1-1}
\Gamma_1(x+1)=x\Gamma_1(x).
\end{equation}
So by Bohr-Mollerup Theorem (see e.g., \cite[p. 44]{GTM172}, the uniqueness of gamma functions), we have
\begin{equation}\label{G1-2}
\Gamma_1(x)=\Gamma(x)R
\end{equation}
for a constant $R$
and by (\ref{gamma}) and (\ref{Riemann2})
\begin{equation}\label{G1-3}
R=\frac{\Gamma_1(1)}{\Gamma(1)}=e^{\zeta'(0,1)}=e^{\zeta'(0)}=e^{-\log\sqrt{2\pi}}=\frac1{\sqrt{2\pi}},
\end{equation}
since $\zeta(s,1)=\zeta(s).$
Therefore we have
\begin{equation}\label{G1-re}
\Gamma_1(x)=\frac{\Gamma(x)}{\sqrt{2\pi}},
\end{equation}
which is what we want.

\section{Euler-Barnes multiple zeta function and the second generalized Lerch formula}\label{EB}

The main aim of this section is to prove the second generalized Lerch's formula (Theorem \ref{gen-lerch2} above).

\subsection*{Proof of Theorem \ref{gen-lerch2}}

The proof goes a similar way as Theorem \ref{gen-lerch}.

The multiple zeta function $\zeta_{E,N}(s,x)$ defined in (\ref{M-Z-2}) may also be represented by the Mellin transform as follows
\begin{equation}\label{B-E-def-3-2}
\zeta_{E,N}(s,x)=\frac1{\Gamma(s)}\int_0^\infty t^se^{-xt}\left(1+e^{-t}\right)^{-N}\frac{dt}{t},
\end{equation}
for Re$(s)>0$ and Re$(x)>0.$

This can be shown as follows.
Start with the  power series 
\begin{equation}\label{re-1-2}
\begin{aligned}
\left(1+e^{-t}\right)^{-N}=\sum_{m_1,\ldots,m_N=0}^\infty (-1)^{m_1+\cdots+m_N} e^{-t(m_1+\cdots+m_N)}, \quad t > 0.
\end{aligned}
\end{equation}
Then substituting  (\ref{re-1-2}) into (\ref{B-E-def-3-2}) we have
\begin{equation}\label{B-E-def-2-2}
\begin{aligned}
\zeta_{E,N}(s,x)
&=\frac1{\Gamma(s)}\int_0^\infty t^se^{-xt}\left(1+e^{-t}\right)^{-N}\frac{dt}{t}\\
&=\sum_{m_1,\ldots,m_N=0}^\infty\frac{(-1)^{m_1+\cdots+m_N}}{\Gamma(s)}\int_0^\infty t^s e^{-(x+m_1+\cdots+m_N) t}\frac{dt}{t} \\
&=\sum_{m_1,\ldots,m_N=0}^\infty\frac{(-1)^{m_1+\cdots+m_N}}{(x+m_1+\cdots+m_N)^s},
\end{aligned}
\end{equation}
for Re$(s)>0$ and Re$(x)>0.$

By (\ref{B-E-def-3-2}) we have
\begin{equation}\label{dir-mul-int-2}
\begin{aligned}
\zeta_{E,N}(s,x)&=\frac1{\Gamma(s)}\int_0^\infty t^se^{-xt}\left(1+e^{-t}\right)^{-N}\frac{dt}{t} \\
&=\frac1{\Gamma(s)}\int_1^\infty\left(\frac1{1+u^{-1}}\right)^{N}u^{-x}(\log u)^{s-1}\frac{du}{u} \\
&=\frac1{\Gamma(s)}\sum_{n=0}^\infty\binom{-N}{n}\int_1^\infty u^{-x-n}(\log u)^{s-1}\frac{du}{u} \\
&\quad(\text{by using (\ref{re-2}) with $w=x+n$}) \\
&=\sum_{n=0}^\infty(-1)^n\binom{n+N-1}{N-1}(x+n)^{-s}.
\end{aligned}
\end{equation}
The above equation is established for $x>0$ and  for all $s\in\mathbb{C}$ by analytic continuation. 
So by (\ref{ru2}) we have
\begin{equation}\label{dir-mul-int-an}
\begin{aligned}
\Gamma_{N}^*(x)=&\exp\left(\frac{\partial}{\partial s}\zeta_{E,N}(s,x)\biggl|_{s=0}\right)\\
&=\exp\left(\sum_{n=0}^\infty(-1)^{n+1}\binom{n+N-1}{N-1}\log(x+n)\right) \\
&=\prod_{n=0}^\infty(x+n)^{(-1)^{n+1}\binom{n+N-1}{N-1}},
\end{aligned}
\end{equation}
which is the desired result.

\subsection*{Proof of Proposition \ref{proposition2}}

Letting $N=1$ in (\ref{M-Z-2}), we recover the alternating Hurwitz zeta function:
$$\zeta_{E}(s,x)=\sum_{n=0}^{\infty}\frac{(-1)^{n}}{(n+x)^{s}}.$$
The following derivative formula of $\zeta_{E}(s,x)$ is shown by Williams and Zhang in \cite[Proposition 3]{WZ} :
$$\zeta_E'\left(0,x\right)=\log\frac{\Gamma\left(\frac{x}{2}\right)}{\Gamma\left(\frac{x+1}{2}\right)}-\frac12\log 2,$$
where $\Gamma(x)$ is the  Euler gamma function.
Hence by (\ref{ru2}) we have
\begin{equation}\label{ex-gamma-2}
\begin{aligned}
\Gamma_{1}^*(x)&=\exp\left(\frac{\partial}{\partial s}\zeta_{E}(s,x)\biggl|_{s=0}\right)\\
&=\frac{\Gamma\left(\frac{x}{2}\right)}{\sqrt{2}\Gamma\left(\frac{x+1}{2}\right)},
\end{aligned}
\end{equation}
which is what we want.


\begin{thebibliography}{99}

\bibitem{Ba}  E.W. Barnes,
\textit{On the theory of the multiple gamma function},
Trans. Cambridge Philos. Soc. \textbf{19} (1904), 374--425.

\bibitem{Choi3} J. Choi, J.R. Quine, \textit{E. W. Barnes' approach of the multiple gamma functions}, J. Korean Math. Soc. \textbf{29} (1992), no. 1, 127--140.


\bibitem{Choi2} J. Choi, 
\textit{Explicit formulas for Bernoulli polynomials of order $n$}, 
Indian J. Pure Appl. Math. \textbf{27} (1996), no. 7, 667--674. 

 

\bibitem{Choi} J. Choi and H.M. Srivastava,
\textit{The multiple Hurwitz zeta function and the multiple Hurwitz-Euler eta function},
Taiwanese J. Math. \textbf{15} (2011), no. 2, 501--522.

\bibitem{Friedmann} T. Friedmann and C.R. Hagen, 
\textit{Quantum mechanical derivation of the Wallis formula for $\pi$}, J. Math. Phys. \textbf{56} (2015), no. 11, 112101, 3 pp. 

\bibitem{arXiv} S. Hu and M.-S. Kim, 
\textit{On $p$-adic multiple Barnes-Euler zeta functions and the corresponding Log Gamma functions}, 
\url{https://arxiv.org/abs/1703.05434}.

\bibitem{KK} N. Kurokawa and S. Koyama,
\textit{Multiple sine functions},
Forum Math. \textbf{15} (2003), no. 6, 839--876.

\bibitem{KW} N. Kurokawa, M. Wakayama, \textit{A generalization of Lerch's formula}, Czechoslovak Math. J. \textbf{54} (129) (2004), no. 4, 941--947. 

\bibitem{KKS} N. Kurokawa, M. Kurihara and T. Saito, 
\textit{Number theory, 3, Iwasawa theory and modular forms}, 
Translated from the Japanese by Masato Kuwata, 
Translations of Mathematical Monographs, 242, Iwanami Series in Modern Mathematics, American Mathematical Society, Providence, RI, 2012. 

\bibitem{Mizuno} Y. Mizuno, \textit{Generalized Lerch formulas: examples of zeta-regularized products}, J. Number Theory \textbf{118} (2006), no. 2, 155--171. 

\bibitem{Miller} S.J. Miller,  
\textit{A probabilistic proof of Wallis's formula for $\pi$}, 
Amer. Math. Monthly 115 (2008), no. 8, 740--745. 

\bibitem{Riemann} B. Riemann, 
\textit{\"Uber die Anzahl der Primzahlen unter einer gegebenen Größe}, 
1859.

\bibitem{Stark} H.M. Stark, \textit{Dirichlet's class-number formula revisited,} A tribute to Emil Grosswald: number theory and related analysis, 571--577, Contemp. Math., 143, Amer. Math. Soc., Providence, RI, 1993. 

\bibitem{GTM172} R. Remmert, 
\textit{Classical topics in complex function theory}, 
Translated from the German by Leslie Kay, Graduate Texts in Mathematics, 172, Springer-Verlag, New York, 1998.

\bibitem{Ru} S.N.M. Ruijsenaars,
\textit{On Barnes' multiple zeta and gamma functions},
Adv. Math.  \textbf{156}  (2000), no. 1, 107--132.
 
 \bibitem{Sondow} J. Sondow,  
 \textit{Analytic continuation of Riemann's zeta function and values at negative integers via Euler's transformation of series}, 
 Proc. Amer. Math. Soc. \textbf{120} (1994), no. 2, 421--424. 
 
 \bibitem{SC} H.M. Srivastava and J. Choi,
\textit{Zeta and $q$-Zeta Functions and Associated Series and Integrals},
Elsevier Science Publishers, Amsterdam, London and New York, 2012.

 
\bibitem{Var} I. Vardi,
\textit{Determinants of Laplacians and multiple Gamma functions},
SIAM J. Math. Anal. \textbf{19} (1988), 493-–507.

\bibitem{Wallis} J. Wallis,  \textit{Arithmetica Infinitorum}, Oxford, England, 1656.

\bibitem{WZ} K.S. Williams and N.Y. Zhang,
\textit{Special values of the Lerch zeta function and the evaluation of certain integrals},
Proc. Amer. Math. Soc. \textbf{119} (1993), no. 1, 35--49.
\end{thebibliography}
\end{document}